\documentclass[reqno]{amsart}
\usepackage{amssymb,url}
\usepackage{hyperref}%,backref
\date{\bf April 11, 2008}
\theoremstyle{plain}
\newtheorem{theorem}{Theorem}
\newtheorem{lemma}{Lemma}
\newtheorem{definition}{Definition}
\newtheorem{corollary}{Corollary}
\newtheorem{proposition}{Proposition}
\theoremstyle{remark}
\newtheorem{remark}{Remark}

\allowdisplaybreaks[4]
\DeclareMathOperator{\T}{\mathbb T}
\DeclareMathOperator{\sigt}{\sigma (t)}
\DeclareMathOperator{\rt}{\rho (t)}

\DeclareMathOperator{\fd}{{\it f}^\Delta (t)}

\DeclareMathOperator{\fs}{{\it f}^\sigma (t)}

\DeclareMathOperator{\Fd}{{\it F}^\Delta (t)}

\begin{document}

\title{An Ostrowski-Gr\"{u}ss type inequality on time scales}

\author[W. Liu]{Wenjun Liu}
\address[W. Liu]{College of Mathematics and Physics\\
Nanjing University of Information Science and Technology \\
Nanjing 210044, China} \email{\href{mailto: W. Liu
<wjliu@nuist.edu.cn>}{wjliu@nuist.edu.cn}}

\author[Q. A. Ng\^{o}]{Qu\^{o}\hspace{-0.5ex}\llap{\raise 1ex\hbox{\'{}}}\hspace{0.5ex}c Anh Ng\^{o}}
\address[Q. A. Ng\^{o}]{Department of Mathematics, Mechanics and Informatics\\
College of Science\\ Vi\d{\^{e}}t Nam National University\\ H\`{a}
N\d{\^{o}}i, Vi\d{\^{e}}t Nam} \email{\href{mailto: Q. A. Ng\^{o}
<bookworm\_vn@yahoo.com>}{bookworm\_vn@yahoo.com}}

\subjclass[2000]{26D15; 39A10; 39A12; 39A13.}

\keywords{Ostrowski inequality, Gr\"{u}ss inequality, Ostrowski-Gr\"{u}ss inequality, time scales.}

\begin{abstract}
In this paper we derive a new inequality of Ostrowski-Gr\"{u}ss type
on time scales and thus unify  corresponding continuous and discrete
versions. We also apply our result to the quantum calculus case.
\end{abstract}

\thanks{This paper was typeset using \AmS-\LaTeX}

\maketitle

\section{introduction}

In 1997, Dragomir and Wang \cite{dw1997} proved the following
Ostrowski-Gr\"{u}ss type integral inequality.
 \begin{theorem}\label{th1}
Let $I\subset \mathbb{R}$ be an open interval, $a, b\in I, a < b.$
If $f: I\rightarrow \mathbb{R}$ is a   differentiable function such
that  there exist constants $\gamma,\Gamma\in \mathbb{R},$ with
 $\gamma\leq f'(x) \leq\Gamma,$ $ x\in [a, b]$. Then we have
\begin{equation} \left|f(x)
-\frac{1}{b-a}\displaystyle\int\limits_a^bf(t)dt-
\frac{f(b)-f(a)}{b-a}\left(x-\frac{a+b}{2}\right)\right|  \leq
\frac{1}{4}(b-a)(\Gamma-\gamma), \label{eq1} \end{equation} for all
$ x\in [a, b]$.
\end{theorem}

This inequality is a connection between the Ostrowski inequality
\cite{mpf} and the Gr\"{u}ss inequality \cite{mpf2}. It can be
applied to bound some special mean and some numerical quadrature
rules. For other related results on the similar integral
inequalities please see the papers \cite{c2001,l1,l3,mpu} and the
references therein.

The aim of this paper is to extend a generalizations of
Ostrowski-Gr\"{u}ss type integral inequality to an arbitrary time
scale.

\section{Time scales essentials}

The development of the theory of time scales was initiated by Hilger
\cite{h1988} in 1988 as a theory capable to contain both difference
and differential calculus in a consistent way. Since then, many
authors have studied the theory of certain integral inequalities on
time scales. For example, we refer the reader to \cite{abp2001,
bm2007, bm2008, OSY, wyy}.

Now we briefly introduce the time scales theory and refer the reader
to Hilger \cite{h1988} and the books \cite{bp2001, bp2003, LSK}  for
further details.

\begin{definition}
{\it A time scale} $\T$ is an arbitrary nonempty closed subset of real numbers.
\end{definition}

%Let $\T$ be a time scale. $\T$ has the topology that it inherits from the real numbers with the standard topology.

\begin{definition}
 For $t \in \T$, we define the {\it forward jump operator} $\sigma : \T \to \T$ by
$
\sigt  = \inf \left\{ {s \in \T:s > t} \right\},
$
while the {\it backward jump operator} $\rho : \T \to \T$ is defined by
$
\rt  = \sup \left\{ {s \in \T:s < t} \right\}.
$
If $\sigt >t $, then we say that $t$ is {\it right-scattered}, while if $\rt < t$ then we say that $t$ is {\it left-scattered}.
\end{definition}

Points that are right-scattered and left-scattered at the same time are called isolated. If $\sigt = t$, the $t$ is called {\it right-dense}, and if $\rt = t$ then $t$ is called {\it left-dense}. Points that are right-dense and left-dense at the same time are called dense.

\begin{definition}
Let $t \in \T$, then two mappings $\mu ,\nu :\T \to \left[ {0, + \infty } \right)$ satisfying
\[
\mu \left( t \right): = \sigt - t, \quad \nu \left( t \right): = t -
\rt
\]
are called the {\it graininess functions}.
\end{definition}

We now introduce the set $\T^\kappa$ which is derived from the time
scales $\T$ as follows. If $\T$ has a left-scattered maximum $t$,
then $\T^\kappa := \T -\{t\}$, otherwise $\T^\kappa := \T$.
Furthermore for a function $f : \T \to \mathbb R$, we define the
function $f^\sigma : \T \to \mathbb R$ by $\fs = f(\sigma(t))$ for
all $t \in \T$.

\begin{definition}
Let $f : \T \to \mathbb R$ be a function on time scales. Then for $t \in \T^\kappa$, we define $\fd$ to be the number, if one exists, such that for all $\varepsilon >0$ there is a neighborhood $U$ of $t$ such that for all $s \in U$
\[
\left| {\fs  - f\left( s \right) - \fd \left( {\sigt - s} \right)} \right| \leqq \varepsilon \left| {\sigt - s} \right|.
\]
We say that $f$ is $\Delta$-differentiable on $\T^\kappa$ provided $\fd$ exists for all $t \in \T^\kappa$.
\end{definition}

%Assume that $f : \T \to \mathbb R$ is a function and let $t \in \T^\kappa$ ($t \ne \min \T$). Then we have the following
%\begin{enumerate}
%  \item[(i)] If $f$ is $\Delta$-differentiable at $t$, then $f$ is continuous at $t$.
%  \item[(ii)] If $f$ is left continuous at $t$ and $t$ is right-scattered, then $f$ is $\Delta$-differentiable at $t$ with
%\[
%\fd = \frac{{\fs - f\left( t \right)}}{{\mu \left( t \right)}}.
%\]
%  \item[(iii)] If $t$ is right-dense, then $f$ is $\Delta$-differentiable at $t$ if and only if
%\[
%\mathop {\lim }\limits_{s \to t} \frac{{f\left( t \right) - f\left( s \right)}}{{t - s}}
%\]
%exists a finite number. In this case
%\[
%\fd = \mathop {\lim }\limits_{s \to t} \frac{{f\left( t \right) - f\left( s \right)}}{{t - s}}
%\]
%  \item[(iv)]If $f$ is $\Delta$-differentiable at $t$, then
%\[
%\fs = f(t) + \mu(t) \fd.
%\]
%\end{enumerate}

\begin{definition} A mapping $f : \T \to \mathbb R$ is called {\it rd-continuous} (denoted by $C_{rd}$) provided if it satisfies
    \begin{enumerate}
      \item $f$ is continuous at each right-dense point or maximal element of $\T$.
      \item The left-sided limit $\mathop {\lim }\limits_{s \to t - } f\left( s \right) = f\left( {t - } \right)$ exists at each left-dense point $t$ of $\T$.
    \end{enumerate}
\end{definition}

\begin{remark}
It follows from Theorem 1.74 of Bohner and Peterson \cite{bp2001} that every rd-continuous function has an anti-derivative.
\end{remark}

\begin{definition} A function $F : \T \to \mathbb R$  is called a $\Delta$-antiderivative
of $f : \T \to \mathbb R$ provided $\Fd =f(t)$ holds for all $t \in \T^\kappa$.
Then the $\Delta$-integral of $f$ is defined by
\[
\int\limits_a^b {f\left( t \right)\Delta t}  = F\left( b \right) - F\left( a \right).
\]
\end{definition}

\begin{proposition} Let $f, g$ be rd-continuous, $a, b, c\in
\mathbb{T}$ and $\alpha, \beta\in \mathbb{R}$. Then
\begin{enumerate}
      \item $
\int\limits_a^b {[\alpha f(t)+\beta g(t)]\Delta t}  =
\alpha\int\limits_a^b {f(t)\Delta t}+\beta\int\limits_a^b
{g(t)\Delta t}, $
      \item $\int\limits_a^b
{f(t)\Delta t}=-\int\limits_b^a {f(t)\Delta t},$
\item $\int\limits_a^b
{f(t)\Delta t}=\int\limits_a^c {f(t)\Delta t}+\int\limits_c^b
{f(t)\Delta t},$
\item $\int\limits_a^b
{f(t)g^\Delta(t)\Delta t}=(fg)(b)-(fg)(a)-\int\limits_a^b
{f^\Delta(t)g(\sigma(t))\Delta t},$
\item $\int\limits_a^a
{f(t)\Delta t}=0.$
 \end{enumerate}
\end{proposition}
\begin{definition}
Let $h_k : \T^2 \to \mathbb R$, $k \in \mathbb N_0$ be defined by
\[
h_0 \left( {t,s} \right) = 1 \quad {\text{ for all }} \quad s,t \in
\T
\]
and then resursively by
\[
h_{k + 1} \left( {t,s} \right) = \int\limits_s^t {h_k \left( {\tau
,s} \right)\Delta \tau }  \quad {\text{ for all }} \quad s,t \in \T.
\]
\end{definition}
%Throughout this paper, we suppose that $\T$ is a time scale, $a, b \in \T$ with $a<b$ and an interval means the intersection of real interval with the given time scale.

The present paper is motivated by the following results: Gr\"{u}ss
inequality on time scales and Ostrowski inequality on time scales
due to Bohner and  Matthews. More precisely, the following so-called
Gr\"{u}ss inequality on time scales was established   in
\cite{bm2007}.

\begin{theorem}[See \cite{bm2007}, Theorem 3.1]\label{th2}
Let $a, b, s\in \T$, $f, g\in C_{rd}$ and $f, g: [a,
b]\rightarrow \mathbb{R}$. Then for
\begin{equation}\label{eq2}
m_1\leq f(s)\leq M_1,\ \ \ \ m_2\leq g(s)\leq M_2,
\end{equation}
we have
\begin{align}
  \left|\frac{1}{b-a}\int\limits_a^b f^\sigma(s)g^\sigma(s)\Delta s-\frac{1}{b-a}\int\limits_a^b f^\sigma(s)\Delta
  s\,
  \frac{1}{b-a}\int\limits_a^b g^\sigma(s)\Delta s\right| \hfill \nonumber\\
   \leq  \frac{1}{4}(M_1-m_1)(M_2-m_2).\label{eq3}
\end{align}
\end{theorem}

The same authors also proved the following so-called Ostrowski
inequality on time scales in \cite{bm2008}.

\begin{theorem}[See \cite{bm2008}, Theorem 3.5]\label{th3}
Let $a, b, s, t\in \T$, $a<b$ and $f: [a, b]\rightarrow
\mathbb{R}$ be differentiable. Then
\begin{equation}
  \left|f(t)-\frac{1}{b-a}\int\limits_a^b f^\sigma(s)\Delta
  s\right|
   \leq  \frac{M}{b-a}(h_2(t,a)+h_2(t,b)),\label{eq4}
\end{equation}
where $M=\sup_{a<t<b}|f^\Delta(t)|.$ This inequality is sharp in the
sense that the right-hand side of (\ref{eq3}) cannot be replaced by
a smaller one.
\end{theorem}

In the present paper we shall first derive  a new inequality of
Ostrowski-Gr\"{u}ss type on time scales by using Theorem \ref{th2}
and then unify corresponding continuous and discrete versions. We
also apply our result to the quantum calculus case.

\section{The Ostrowski-Gr\"{u}ss type inequality on time scales}

Similarly as in \cite{dw1997}, the
Ostrowski-Gr\"{u}ss type inequality can be shown for
general time scales.

\begin{theorem}\label{th4}
Let $a, b, s, t\in \T$, $a<b$ and $f: [a, b]\rightarrow
\mathbb{R}$ be differentiable. If $f^\Delta$ is rd-continuous and
$$\gamma\leq f^\Delta(t)\leq \Gamma,\ \ \ \forall\ t\in [a, b].$$ Then we have
\begin{equation}
  \left|f(t)-\frac{1}{b-a}\int\limits_a^b f^\sigma(s)\Delta
  s-
\frac{f(b)-f(a)}{(b-a)^2}\left( h_2(t,a)- h_2 (t, b) \right)\right|
   \leq  \frac{1}{4}(b-a)(\Gamma-\gamma),\label{eq5}
\end{equation}
for all $t\in [a, b]$.
\end{theorem}

To prove Theorem \ref{th4}, we need the following Montgomery
Identity.
\begin{lemma}[Montgomery Identity, see \cite{bm2008}]\label{le1}
Let $a, b, s, t\in \T$, $a<b$ and $f: [a, b]\rightarrow
\mathbb{R}$ be differentiable. Then
\begin{equation}
 f(t)=\frac{1}{b-a}\int\limits_a^b f^\sigma(s)\Delta
  s+\frac{1}{b-a}\int\limits_a^b p(t,s)f^\Delta(s)\Delta
  s,\label{eq6}
\end{equation}
where
\begin{equation}
 p(t,s)=\left\{ \begin{array}{ll}
  s-a,\ \ a\leq s<t, \hfill \nonumber\\
   s-b,\ \ t\leq s\leq b.
\end{array} \right.\label{eq7}
\end{equation}
\end{lemma}

\begin{proof}[Proof of Theorem~\ref{th4}]
By applying Lemma \ref{le1}, we get
\begin{equation}
 f(t)-\frac{1}{b-a}\int\limits_a^b f^\sigma(s)\Delta
  s=\frac{1}{b-a}\int\limits_a^b p(t,s)f^\Delta(s)\Delta
  s,\label{eq7}
\end{equation}
for all $t\in [a, b]$. It is clear that for all   $t\in [a, b]$ and
$s\in [a, b]$ we have
$$t-b\leq p(t,s)\leq t-a.$$
Applying Theorem \ref{th2} to the  mapping $p(t, \cdot)$ and
$f^\Delta(\cdot)$, we get
\begin{align}
  \left|\frac{1}{b-a}\int\limits_a^b p(t,s)f^\Delta(s)\Delta
  s\right.&\left.-\frac{1}{b-a}\int\limits_a^b p(t,s) \Delta s\,
  \frac{1}{b-a}\int\limits_a^b f^\Delta(s)\Delta s\right|   \nonumber\\
  & \leq  \frac{1}{4}[(t-a)-(t-b)](\Gamma-\gamma)  \nonumber\\
  & \leq \frac{1}{4}(b-a)(\Gamma-\gamma).\label{eq8}
\end{align}
By a simple calculation we get
\begin{align*}\int\limits_a^b p(t,s) \Delta s&=\int\limits_a^t (s-a) \Delta s+\int\limits_t^b (s-b) \Delta
s\\
&=\int\limits_a^t (s-a) \Delta s-\int\limits_b^t (s-b) \Delta s \\&=
h_2(t,a) - h_2 (t, b)\end{align*} and
$$ \frac{1}{b-a}\int\limits_a^b f^\Delta(s)\Delta s=\frac{f(b)-f(a)}{b-a}.$$
By combining (\ref{eq7}), (\ref{eq8}) and the above two equalities,
we obtain (\ref{eq5}).
\end{proof}

If we apply the Ostrowski-Gr\"{u}ss type inequality
to different time scales, we will get some well-known and some new
results.

\begin{corollary}\label{co1} {\rm(Continuous case)}. Let $\mathbb{T }= \mathbb{R}$.
Then our delta integral is the usual Riemann integral from calculus. Hence,
\[
h_2 \left( {t,s} \right) = \frac{{\left( {t - s} \right)^2 }}{2}, \quad {\text{ for all }} \quad t, s \in \mathbb R.
\]
This leads us to state the following inequality
\begin{equation} \left|f(t)
-\frac{1}{b-a}\displaystyle\int\limits_a^bf(s)ds-
\frac{f(b)-f(a)}{b-a}\left(t-\frac{a+b}{2}\right)\right|  \leq
\frac{1}{4}(b-a)(\Gamma-\gamma), \label{eq9} \end{equation}  for all
$t\in [a, b]$,  where $\gamma\leq f'(t) \leq\Gamma,$ which is
exactly the Ostrowski-Gr\"{u}ss type inequality
shown in Theorem \ref{th1}.
\end{corollary}

\begin{corollary}\label{co2} {\rm (Discrete case)}. Let $\mathbb{T }= \mathbb{Z}$, $a=0$, $b=n$, $s=j$, $t=i$ and $f(k)=x_k$. With these, it is known that
\[
h_k \left( {t,s} \right) = \left( {\begin{array}{*{20}c}
   {t - s}  \\
   k  \\
 \end{array} } \right) , \quad {\text{ for all }} \quad  t,s \in \mathbb Z.
\]
Therefore,
\[
h_2 \left( {t,0} \right) = \left( {\begin{array}{*{20}c}
   t  \\
   2  \\
 \end{array} } \right) = \frac{{t\left( {t - 1} \right)}}
{2} \quad , \quad h_2 \left( {t,n} \right) = \left( {\begin{array}{*{20}c}
   {t - n}  \\
   2  \\
 \end{array} } \right) = \frac{{\left( {t - n} \right)\left( {t - n - 1} \right)}}
{2}.
\]
Thus, we have
\begin{equation} \left|x_i
-\frac{1}{n}\sum_{j=1}^nx_j-
\frac{x_n}{n}\left(i-\frac{n+1}{2}\right)\right|  \leq
\frac{1}{4}n(\Gamma-\gamma), \label{eq10} \end{equation}  for all
 $i=\overline{1,n}$,  where $\gamma\leq \Delta x_i \leq\Gamma$.
\end{corollary}

\begin{corollary}\label{co3} {\rm (Quantum calculus case)}. Let
$\mathbb{T }= q^{\mathbb{N}_0}$, $q>1$, $a=q^m, b=q^n$ with  $m<n$. In this situation, one has
\[
h_k \left( {t,s} \right) = \prod\limits_{\nu  = 0}^{k - 1} {\frac{{t - q^\nu  s}}
{{\sum\limits_{\mu  = 0}^\nu  {q^\mu  } }}}, \quad {\text{ for all }} \quad  t,s \in \T.
\]
Therefore,
\[
h_2 \left( {t,q^m } \right) = \frac{{\left( {t - q^m } \right)\left( {t - q^{m + 1} } \right)}}
{{1 + q}} \quad , \quad h_2 \left( {t,q^n } \right) = \frac{{\left( {t - q^n } \right)\left( {t - q^{n + 1} } \right)}}
{{1 + q}}.
\]
Then
\begin{equation}\label{eq11}
\begin{split}
  \Bigg|f(t)-\frac{1}{q^n-q^m}\int\limits_{q^m}^{q^n} f^\sigma(s)\Delta  s &-\frac{f(q^n)-f(q^m)}{q^n-q^m}\left(t-\frac{q^{2n+1} - q^{2m+1}}{q+1}\right)\Bigg| \\
&\leq  \frac{1}{4}(q^n-q^m)(\Gamma-\gamma),
\end{split}
\end{equation}
where $$\gamma\leq  \frac{f(qt)-f(t)}{(q-1)(t)} \leq\Gamma,\ \ \forall \ t\in [a,
b].$$
\end{corollary}

If $f^\Delta$ is bounded on $[a, b]$ then we have the following
corollary.
\begin{corollary}\label{co4}  With the assumptions in Theorem~\ref{th4},
if $|f^\Delta(t)|\leq M$ for all $t\in [a, b]$ and some
positive constant $M$, then we have
\begin{equation}
 \left|f(t)-\frac{1}{b-a}\int\limits_a^b f^\sigma(s)\Delta
  s-
\frac{f(b)-f(a)}{(b-a)^2}\left( h_2(t,a)- h_2 (t, b) \right)\right|
   \leq  \frac{1}{2}(b-a)M,\label{eq12}
\end{equation}
for all $t\in [a, b]$.
\end{corollary}

Furthermore, choosing $t = (a + b)/2$ and $t = b$, respectively, in
(\ref{eq5}), we have the following corollary.

\begin{corollary}\label{co5}  With the assumptions in Theorem~\ref{th4},  we have
\begin{align}
  &\ \ \ \left|f\left(\frac{a+b}{2}\right)-\frac{f(b)-f(a)}{(b-a)^2}\left[
  h_2\left(\frac{a+b}{2},a\right)- h_2\left(\frac{a+b}{2},b\right)\right]
  -\frac{1}{b-a}\int\limits_a^b f^\sigma(s)\Delta
  s\right|   \nonumber\\
  &  \leq  \frac{1}{4}(b-a)(\Gamma-\gamma)
  \ \ \ \ \ \ \ \ \ \ \ \ \ \left(\mbox{if}\ \ \frac{a+b}{2}\in \mathbb{T}\right)
  \label{eq13}
\end{align}
and
\begin{equation}
  \left|f(b)-\frac{f(b)-f(a)}{(b-a)^2}h_2(b,a)-\frac{1}{b-a}\int\limits_a^b f^\sigma(s)\Delta
  s\right|
   \leq  \frac{1}{4}(b-a)(\Gamma-\gamma).\label{eq14}
\end{equation}
\end{corollary}

\medskip
\section*{Acknowledgements}

This work was supported by the Science Research Foundation of
Nanjing University of Information Science and Technology and the
Natural Science Foundation of Jiangsu Province Education Department
under Grant No.07KJD510133.

\end{document}